\documentclass[11pt]{amsart}

\usepackage{amssymb,amsmath,amsthm,amscd,mathrsfs,graphicx}
\usepackage{amscd}
\usepackage{amsxtra}
\usepackage[cmtip,all]{xy}
\usepackage{pb-diagram,pb-xy}

\begin{document}

\theoremstyle{plain} \numberwithin{equation}{section}
\newtheorem*{convention}{}
\newtheorem*{unlabeledtheorem}{Theorem}
\newtheorem{theorem}{Theorem}[section]
\newtheorem{question}{Question}[section]
\newtheorem{theoremalpha}{Theorem}
\renewcommand{\thetheoremalpha}{\Alph{theoremalpha}}
\newtheorem{corollary}[theorem]{Corollary}
\newtheorem{lemma}[theorem]{Lemma}
\newtheorem{proposition}[theorem]{Proposition}
\newtheorem{conjecture}{Conjecture}

\theoremstyle{definition}
\newtheorem{definition}[theorem]{Definition}

\theoremstyle{remark}
\newtheorem{hypothesis}{Hypothesis}
\newtheorem{remark}[theorem]{Remark}
\newtheorem*{unlabeledremark}{Remark}
\newtheorem{example}[theorem]{Example}

\newcommand{\Spec}{\text{Spec}}

\newcommand{\dual}{\omega} 

\newcommand{\isoto}{\stackrel{\backsim}{\longrightarrow}}

\newcommand{\nSm}{\operatorname{D}} 

\newcommand{\ShAut}{\underline{\operatorname{Aut}}}
\newcommand{\Aut}{\operatorname{Aut}}
\newcommand{\ShHom}{\underline{\operatorname{Hom}}}
\newcommand{\Hom}{\operatorname{Hom}}
\newcommand{\ShEnd}{\underline{\operatorname{End}}}
\newcommand{\End}{\operatorname{End}}
\newcommand{\Ext}{\operatorname{Ext}}
\newcommand{\ShExt}{\underline{\operatorname{Ext}}}

\newcommand{\Quot}{\text{Quot}_{\dual}}
\newcommand{\Hilb}{\text{Hilb}_{X}}
\newcommand{\Jac}{\overline{J}_{X}}

\newcommand{\AbelH}{\operatorname{A}_{\text{h}}}
\newcommand{\AbelQ}{\operatorname{A}_{\text{q}}}

\newcommand{\bbA}{\mathbf{A}}
\newcommand{\bbB}{\mathbf{B}}
\newcommand{\bbC}{\mathbf{C}}
\newcommand{\bbD}{\mathbf{D}}
\newcommand{\bbE}{\mathbf{E}}
\newcommand{\bbF}{\mathbf{F}}
\newcommand{\bbG}{\mathbf{G}}
\newcommand{\bbH}{\mathbf{H}}
\newcommand{\bbI}{\mathbf{I}}
\newcommand{\bbJ}{\mathbf{J}}
\newcommand{\bbK}{\mathbf{K}}
\newcommand{\bbL}{\mathbf{L}}
\newcommand{\bbM}{\mathbf{M}}
\newcommand{\bbN}{\mathbf{N}}
\newcommand{\bbO}{\mathbf{O}}
\newcommand{\bbP}{\mathbf{P}}
\newcommand{\bbQ}{\mathbf{Q}}
\newcommand{\bbR}{\mathbf{R}}
\newcommand{\bbS}{\mathbf{S}}
\newcommand{\bbT}{\mathbf{T}}
\newcommand{\bbU}{\mathbf{U}}
\newcommand{\bbV}{\mathbf{V}}
\newcommand{\bbW}{\mathbf{W}}
\newcommand{\bbX}{\mathbf{X}}
\newcommand{\bbY}{\mathbf{Y}}
\newcommand{\bbZ}{\mathbf{Z}}

\newcommand{\calA}{{ \mathcal A}}
\newcommand{\calB}{{ \mathcal B}}
\newcommand{\calC}{{ \mathcal C}}
\newcommand{\calD}{{ \mathcal D}}
\newcommand{\calE}{{ \mathcal E}}
\newcommand{\calF}{{ \mathcal F}}
\newcommand{\calG}{{ \mathcal G}}
\newcommand{\calH}{{ \mathcal H}}
\newcommand{\calI}{{ \mathcal I}}
\newcommand{\calJ}{{ \mathcal J}}
\newcommand{\calK}{{ \mathcal K}}
\newcommand{\calL}{{ \mathcal L}}
\newcommand{\calM}{{ \mathcal M}}
\newcommand{\calN}{{ \mathcal N}}
\newcommand{\calO}{{ \mathcal O}}
\newcommand{\calP}{ {\mathcal P} }
\newcommand{\calQ}{{ \mathcal Q}}
\newcommand{\calR}{{ \mathcal R}}
\newcommand{\calS}{{ \mathcal S}}
\newcommand{\calT}{{ \mathcal T}}
\newcommand{\calU}{{ \mathcal U}}
\newcommand{\calV}{{ \mathcal V}}
\newcommand{\calW}{{ \mathcal W}}
\newcommand{\calX}{{\mathcal X}}
\newcommand{\calY}{{\mathcal Y}}
\newcommand{\calZ}{ {\mathcal Z}}

\makeatletter
\@namedef{subjclassname@2010}{ \textup{2010} Mathematics Subject 
Classification} \makeatother

\title[Non-Smoothable Component]{An Explicit Non-Smoothable Component of the Compactified Jacobian}
\author[Kass]{Jesse Leo Kass}
\address{University of Michigan, Department of Mathematics,  530 Church St, Ann Arbor, MI 48103, USA}
\email{kass@math.harvard.edu}

\date{\today}

\thanks{J.  L.  K. was supported by NSF grant DMS-0502170.}
\date{\today}

\bibliographystyle{amsalpha}

\begin{abstract}
	This paper studies the components of the moduli space of rank $1$, torsion-free sheaves, or compactified Jacobian, of a non-Gorenstein curve.  We exhibit a generically reduced component of dimension equal to the arithmetic genus and prove that it is the only non-smoothable component when the curve has a unique singularity that is of finite representation type.  Analogous results are proven for the Hilbert scheme of points and the Quot scheme parameterizing quotients of the dualizing sheaf.
\end{abstract}

\maketitle

\section{Introduction}
Associated to a curve $X$ is the \texttt{compactified Jacobian} $\Jac^{d}$, or moduli space of rank $1$, torsion-free sheaves of degree $d$.  Examples of rank $1$, torsion-free sheaves are line bundles, and the closure of the corresponding line bundle locus in $\Jac^{d}$ is an irreducible component called the \texttt{smoothable} component.  

Are there other components?  Altman--Iarrobino--Kleiman \cite{altman77} and Kleiman--Kleppe \cite{kleppe81a} have answered this question.  They showed that there are other components (i.e.~\texttt{non-smoothable} components) precisely when $X$ has a non-planar singularity.  This paper is concerned with the natural follow-up question: when $X$ has a non-planar singularity, what are the additional components? We prove two theorems addressing this question.

If $X$ is a curve with a non-Gorenstein singularity $p_0 \in X$ and $d \in \bbZ$ is an integer, then we set $\nSm_d \subset \Jac^{d}$ equal to the subset  corresponding to sheaves that become isomorphic to the dualizing sheaf $\dual$  upon passing to the completed local ring $\widehat{\calO}_{X,p_0}$. (See Def.~\ref{Def: NonSm}.)

\begin{theoremalpha} \label{Thm: CompsExist}
	Let $X$ be a curve with a unique non-Gorenstein singularity $p_0 \in X$.  Then the closure of $\nSm_{d}$ is a non-smoothable component of $\Jac^{d}$. 
\end{theoremalpha}
This is a special case of Theorem~\ref{Thm: CompsExistRefined}.  The latter includes  analogous results about the \texttt{Hilbert scheme} $\Hilb^{d}$ and the \texttt{Quot scheme} $\Quot^{d}$. In the above statement, we only assume $X$ has a unique non-Gorenstein singularity to simplify exposition.  When $X$ has $n$ non-Gorenstein singularities, a modification of Theorem~\ref{Thm: CompsExist} produces $2^n-1$ non-smoothable components.  

 Theorem~\ref{Thm: CompsExist} accounts for all the non-smoothable components when $X$ has a unique singularity that is of finite representation type.
\begin{theoremalpha} \label{Thm: AllComps}
	Let $X$ be a curve with a unique singularity that is of finite representation type and non-Gorenstein.  Then $\Jac^{d}$ has exactly two components: the smoothable component and the closure of $\nSm_{d}$.
\end{theoremalpha}
As in Theorem~\ref{Thm: CompsExist}, we only assume $X$ has a unique singularity in order to simplify the exposition.  Theorem~\ref{Thm: AllComps} is deduced from Theorem~\ref{Thm: RefinedAllComps}, and the latter also establishes the analogous result for $\Quot^{d}$.  

Recall that a singularity $p_0 \in X$ is said to be of \texttt{finite representation type} if there are only finitely many isomorphism classes of maximal Cohen-Macaulay modules over
$\widehat{\calO}_{X,p}$.  Finite representation type is a strong condition to impose.  The singularities that are of finite representation type and planar  are exactly the ADE curve singularities \cite{greuel85}.   The non-planar curve singularities of finite representation type are all  non-Gorenstein and fall into one infinite family and three exceptional cases.  These singularities are listed in Table~\ref{Table: CurveSing} (p.~\pageref{Table: CurveSing}), and their classification is  discussed in Section~\ref{Sec: Tables}. 

The rank $1$, torsion-free modules over a singularity that is of finite representation type are listed in Table~\ref{Table: BigTable} (p.~\pageref{Table: BigTable}), and their classification can be summarized succinctly.  Given the ring $\calO$ of a singularity from the table and a finite extension $\calO' \supset \calO$ contained in  $\operatorname{Frac} \calO$, both the over-ring $\calO'$ and its dualizing module $\dual'$ are rank $1$, torsion-free $\calO$-modules.  Classification shows that if $I$ is a rank $1$, torsion-free $\calO$-module, then there exists a unique  extension $\calO' \supset \calO$ such that  $I$ is isomorphic  to either $\calO'$ or $\dual'$.

\subsection*{Comparison with Past Work}
The author believes this paper provides the first complete enumeration of the irreducible components of a reducible compactified Jacobian.  The proof that $\Jac^d$ is reducible for $X$ non-planar was given by Kleiman--Kleppe in \cite{kleppe81a}.  Let $g$ be the arithmetic genus of $X$ and $e$ the minimal number of generators of the stalk of $\dual$ at $p_0 \in X$.  The authors of \cite{kleppe81a} exhibit a $(g+e-2)$-dimensional locus in $\Jac^d$ with the property that the general element is not a line bundle \cite[Prop.~4]{kleppe81a}. When $e \ge 2$, dimensional considerations show that this locus must be contained in some non-smoothable component.   When $X$ has a unique singularity that is of finite representation type and non-Gorenstein (hence $e=2$),
the locus constructed by Kleiman--Kleppe coincides with the closure of $\nSm_d$.  However, a comparison of dimensions shows that the closure of $\nSm_d$ cannot contain the Kleiman--Kleppe locus when $e \ge 3$.  For more general surveys of compactified Jacobians, Hilbert schemes, and related topics, the author directs the reader to \cite{iarrobino87} and \cite{cartwright}.

\subsection*{Conventions}  We work over a fixed algebraically closed field $k$ of characteristic $0$.  If $t \in T$ is a point of a $k$-scheme, then we write $k(t)$ for the residue field of the local ring $\calO_{T,t}$ and call this field the \textbf{fiber}.  A \textbf{curve} is an integral projective $k$-scheme of dimension $1$. The \textbf{local ring of a curve singularity} is the completed local ring $\widehat{\calO}_{X,p}$ of a curve at some closed point $p$.

	Given a coherent $\calO_{X}$-module $F$ on a $k$-scheme $X$, we write $F_{x}$ for the \textbf{stalk} of $F$ at $x$.  We call the $k(x)$-module $F_{x} \otimes_{\calO_{X,x}} k(x)$ the \textbf{fiber}.   If $X$ is a curve, then we define the \textbf{degree} of $F$ by $\chi(X,F) = \deg(F) + \chi(X, \calO_{X})$.  We say that $F$ is a \textbf{maximal Cohen-Macaulay} (or \textbf{maximal CM}) sheaf if the localization $F_p$ has depth $1$ for all closed points $p \in X$.  If $F$ additionally has the property that the generic rank is $1$, then we say that $F$ is a \textbf{rank $1$, torsion-free sheaf}.   We write $\dual$ for the \textbf{dualizing sheaf} of $X$, which is a rank $1$, torsion-free sheaf.

	Given two coherent $\calO_{X}$-modules $F$ and $G$, the $\calO_{X}$-module whose sections over an open $U \subset X$ are homomorphisms $F|_{U} \to G|_{U}$ is denoted by  $\ShHom(F,G)$.  Set $F^{\vee} := \ShHom(F, \calO_{X})$.   	
	
	The \textbf{compactified Jacobian} $\Jac^{d}$ of degree $d$ is the projective $k$-scheme that parameterizes rank $1$, torsion-free sheaves of degree $d$ on $X$.  (See \cite[Def.~5.11, Thm.~8.1]{altman80} for the precise definition.)  We write $[I] \in \Jac^{d}$ for the closed point corresponding to a rank $1$, torsion-free sheaf $I$.
	
	The \textbf{Quot scheme} $\Quot^{d}$ of degree $d$ is the projective $k$-scheme that parameterizes rank $d$ quotients $q \colon \dual \twoheadrightarrow Q$ of the dualizing sheaf $\dual$.  We write $[q]$ for the closed point corresponding to $q$.  Similarly, the \textbf{Hilbert scheme} $\Hilb^{d}$ of degree $d$ is the projective $k$-scheme that parameterizes rank $d$ quotients of the structure sheaf $\calO_{X}$ or, equivalently, degree $d$ closed subschemes $Z \subset X$ of $X$.  We write $[Z] \in \Hilb^{d}$ for the closed point corresponding to $Z$.  (Precise definitions can be found in \cite[Def.~2.5, Thm.~2.6]{altman80}.)  Both the kernel $\ker(q)$ associated to $[q] \in \Quot^{d}$ and the ideal $I_{Z}$ associated to $[Z] \in \Hilb^{d}$ are rank $1$, torsion-free sheaves.


\section{Proof of Theorem~\ref{Thm: CompsExist}} \label{Sec: PfOfThmA}
Here we prove Theorem~\ref{Thm: CompsExist}.  We begin by recording the definition of the non-smoothable locus $\nSm_{d}$ more formally.
\begin{definition} \label{Def: NonSm}
	If $x \in \Jac^{d}$ is a (possibly non-closed) point, then set $\overline{k}(x)$ equal to an algebraic closure of  the residue field of $\Jac^{d}$ at $x$. 
	
	Set $I_{x}$ equal to the rank $1$, torsion-free sheaf on $X \otimes \overline{k}(x)$ that is the pullback of a universal family on $X \times \Jac^{d}$ under $X \otimes \overline{k}(x) \to X \times \Jac^{d}$.
	
	We define $$\nSm_{d} \subset \Jac^{d}$$ to be the subset of points $x \in \Jac^{d}$ such that $I_{x}$ and the dualizing sheaf $\dual \otimes \overline{k}(x)$ become isomorphic after tensoring with $\widehat{\calO}_{X \otimes \overline{k}(x), p_{0}}$.
\end{definition}

In addition to proving that $\nSm_{d} \subset \Jac^{d}$ is non-smoothable, we also prove similar results  about the Hilbert scheme $\Hilb^{d}$ parameterizing closed subschemes and the Quot scheme $\Quot^{d}$ parameterizing quotients of the dualizing sheaf $\dual$ (Thm.~\ref{Thm: CompsExistRefined}).  The schemes $\Hilb^{d}$ and $\Quot^{d}$ are related to $\Jac^{d}$ by an Abel map, the properties of which were studied in \cite{altman80}.  

The Abel map $\AbelQ \colon \Quot^{d} \to \Jac^{(2g-2)-d}$ is defined by the rule 
\begin{displaymath}
	[q] \in \Quot^{d} \mapsto [\ker(q)] \in \Jac^{(2g-2)-d}.
\end{displaymath}	
This map fibers $\Quot^{d}$ by projective spaces of possibly varying dimension: the fiber over $[I] \in \Jac^{(2g-2)-d}$ is $\bbP \Hom(I, \dual)$.  This projective space is non-empty once $d \ge g$ and of dimension $d-g$ when $d \ge 2g -1$ \cite[Thm.~8.4]{altman80}.  

The Hilbert scheme $\Hilb^{d}$ also admits an Abel map.  The rule $[Z] \mapsto [I_{Z}]$ defines a morphism  $\AbelH \colon \Hilb^{d} \to \Jac^{-d}$ with the property that the fiber over $[I]$ is the projective space $\bbP \Hom(I, \calO_{X})$.  One important difference between $\AbelQ$ and $\AbelH$ is that the fibers of $\AbelH$ may not  all be of dimension $d-g$ once $d \ge 2g-1$.  Indeed, when $X$ is Gorenstein, $\Hilb^{d}$ can be identified with $\Quot^{d}$ in a way that respects Abel maps, so this condition on the fibers does hold, but in the non-Gorenstein case, the dimension of a fiber of $\AbelH$ is non-constant as a function of the base.  This problem does not arise if we restrict our attention to $\nSm_{-d} \subset \Jac^{-d}$ as will be shown by the following series of lemmas.

\begin{lemma} \label{Lemma: LocClosed}
		Let $X$ be a genus $g$ curve with a unique non-Gorenstein singularity $p_0$.  Then $\nSm_d \subset \Jac^{d}$ is a $g$-dimensional, irreducible locally closed subset.
\end{lemma}
\begin{proof}
	This lemma can be deduced from results in \cite{gagne00}, but we include a proof for the sake of completeness.  We begin by showing that  $\nSm_{d}$ is the image of $J^{d-(2g-2)}_{X}$ under the map $L \mapsto \dual \otimes L$ and that this map is injective.  Certainly the image is contained in $\nSm_{d}$.  To  establish the reverse inclusion, we must show that, given $x \in \nSm_{d}$, the sheaf $I_{x}$ on $X \otimes \overline{k}(x)$ is isomorphic to $\left(\dual \otimes \overline{k}(x)\right) \otimes L$ for a  line bundle $L$.	
	To construct $L$, observe that an examination of completed stalks shows that
	\begin{displaymath}
		M := \ShHom(I_{x},\dual\otimes \overline{k}(x))
	\end{displaymath}
	is a line bundle.  However, $I_{x}$ is $\dual$-reflexive (by \cite[2.2.1]{gagne00}), so 
	\begin{align*}
		I_{x}=& \ShHom(M,\dual \otimes \overline{k}(x)) \\
		=& \left(\dual \otimes \overline{k}(x)\right) \otimes M^{-1},
	\end{align*}
	and  we take $L=M^{-1}$.  This proves the reverse inclusion.  Furthermore, the construction shows that $L$ is unique, so $L \mapsto L \otimes \dual$ is injective.   (We thank the anonymous referee for suggesting this argument.)

	A restatement of this description of $\nSm_{d}$ is: if we fix base points, then we can identify $\nSm_{d}$ with the orbit of a point under a fixed-point free action of the Jacobian $J_{X}^{0}$.  Because $X$ is a curve,  $J^{0}_{X}$ is a reduced and irreducible group scheme \cite[5.23]{kleiman05}, and it is a result of Chevalley that the orbit of any such group scheme is locally closed \cite[2.3.3]{springer}.  We can conclude that $\nSm_{d}$ is an irreducible locally closed subset of dimension $ \dim J^{0}_{X} =g$.  This completes the proof.
\end{proof}
Because $\nSm_{d} \subset \Jac^{d}$ is a locally closed subset, it has a natural subscheme structure --- the reduced subscheme structure.  For the remainder of this article, we will consider $\nSm_{d}$ as a subscheme rather than  as a subset.   The proof of Lemma~\ref{Lemma: LocClosed} shows that $\nSm_{d}$ is isomorphic to $J^{0}_{X}$ as a scheme.

\begin{definition}
	Given $d \in \bbZ$, define integers
		\begin{align*}
		d_0 	:=& 1-\deg(\dual^{\vee}), \\
		n_{d} 	:=& d+2g-2+\deg(\dual^{\vee})\\
			=& d-d_0+2g-1.
	\end{align*}
\end{definition}

\begin{lemma}
	If $d \ge d_0-g+1$ (i.e.~$n_{d}-g \ge 0$), then $\dim \Hom(I,\calO_{X}) \ge n_{d}-g+1$ for all $[I] \in \nSm_{-d}$.  Furthermore, equality holds for all $I$ when $d \ge d_0$ and for some $I$ when $d$ is arbitrary.
\end{lemma}
\begin{proof}
	Suppose first $d \ge d_{0}$.   Given $[I] \in \nSm_{-d}$, the proof of the previous lemma shows that we can write 
	$I = L \otimes \dual$ for some line bundle $L$ of  degree $-d-(2g-2)$.  Then
\begin{align*}
	\dim \Hom(I, \calO_{X}) 	=& 	\dim \Hom( \dual \otimes L, \calO_{X}) \\
									=& \dim H^{0}( \ShHom( \dual \otimes L, \calO_{X})) \\
									=& \chi( \ShHom( \dual \otimes L, \calO_{X})) + \dim H^{1}( \ShHom(\dual \otimes L, \calO_{X})).
\end{align*}
The degree of  $\ShHom( \dual \otimes L,\calO_{X}) = \dual^{\vee} \otimes L^{-1}$ is strictly larger than $2g-2$, so this sheaf has no higher cohomology (by \cite[Prop.~3.5(iii)(g)]{altman80}).  Elementary algebra shows that $\dim \Hom(I,\calO_{X}) = n_{d}-g+1$.

	Now suppose $d_0 > d \ge d_0-g+1$.  Then every element $[I] \in \nSm_{-d}$ can be written as 
$I = J \otimes L$ for $L$ a line bundle of degree $d_0-d$ and $[J] \in \nSm_{-d_0}$.  Computing as before, we have
\begin{align*}
	\dim \Hom(I, \calO_{X}) =& \chi(J^{\vee})-\deg(L)+ \dim H^{1}(J^{\vee} \otimes L^{-1}) \\
				\ge& \chi(J^{\vee})-\deg(L) \\
				=&  n_{d}-g+1.
\end{align*}
This establishes the desired lower bound. A sheaf achieving this lower bound is $J(p_1 + \dots + p_{d_0-d})$ for
 $[J] \in \nSm_{-d_0}$ arbitrary and $p_1, \dots, p_{d_0-d} \in X$ general (as vanishing at a general point is a non-trivial linear condition).  
\end{proof}

\begin{definition}
	If $d \ge d_0-g+1$,  set
	\begin{displaymath}
		\nSm_{-d}^{\operatorname{o}} := \{ [I] \in \nSm_{-d} \colon \dim \Hom(I,\calO_{X}) = n_{d}-g+1 \}.
	\end{displaymath}
\end{definition}

\begin{lemma} \label{Lemma: AbelForHilb}
	Let $X$ be a genus $g$ curve with a unique non-Gorenstein singularity $p_0$.  If $d \ge d_0-g+1$, then $\nSm^{\operatorname{o}}_{-d} \subset \nSm_{-d}$ is open and the restriction of the Abel map
	\begin{equation} \label{Eqn: RestrictAbel}
		\AbelH \colon \AbelH^{-1}(\nSm^{\operatorname{o}}_{-d}) \to \nSm^{\operatorname{o}}_{-d}
	\end{equation}
	is smooth with fibers isomorphic to $\bbP^{n_{d}-g}$. 
\end{lemma}
\begin{proof}
	We prove this by using the description of $\Hilb^{d}$ as the projectivization of a coherent sheaf on $\Jac^{-d}$.  Recall that if we choose a universal family $I_{\text{uni}}$ on $X \times \Jac^{-d}$, then $\Hilb^{d} = \bbP(\operatorname{H})$ for $\operatorname{H} = \operatorname{H}(I_{\text{uni}}, \calO_{X \times \Jac^{-d}})$ (defined in \cite[Sec.~1]{altman80}).  Given $y \in \nSm_{-d}$, the fiber $A^{-1}(y) = \bbP \Hom(I, \calO_{X})$ equals $\bbP \Hom( \operatorname{H}\otimes k(y), k)$ by \cite[1.1.1]{altman80}.   This fiber is thus isomorphic to $\bbP^{n_{d}-g}$.  precisely when $\dim \Hom(I,\calO_{X})=n_{d}-g+1$ or, in other words, $y \in \nSm_{-d}^{\operatorname{o}}$.  Furthermore, $n_{d}-g+1$ is the minimal possible value of  $\dim_{k} \operatorname{H}\otimes k(y)$ for $y \in \nSm_{-d}$ by the previous lemma.  We can conclude that $\nSm_{-d}^{\operatorname{o}} \subset \nSm_{-d}$ is open.
	
	To complete the proof, we must show that the restriction of $\AbelH$ to $\AbelH^{-1}(\nSm^{\operatorname{o}}_{-d})$ is smooth.  The scheme $\nSm_{-d}^{\operatorname{o}}$ is reduced as it inherits this property from $\nSm_{-d}$ (which was defined to be reduced).  Now consider the sheaf $\operatorname{H} \otimes \calO_{\nSm_{-d}^{\operatorname{o}}}$ on $\calO_{\nSm_{-d}^{\operatorname{o}}}$.  This sheaf has the property that all the fibers  $(\operatorname{H} \otimes \calO_{\nSm_{-d}^{\operatorname{o}}}) \otimes k(y)$ have the same dimension.  Because $\nSm_{-d}^{\operatorname{o}}$ is reduced, we can conclude that $\operatorname{H} \otimes \calO_{\nSm_{-d}^{\operatorname{o}}}$ is locally free.  In particular, $\AbelH \colon \AbelH^{-1}(\nSm^{\operatorname{o}}_{-d}) \to \nSm^{\operatorname{o}}_{-d}$ is the projectivization of a locally free sheaf, hence smooth.  This completes the proof.
\end{proof}

We now deduce Theorem~\ref{Thm: CompsExist} in its most general form.
\begin{theorem} \label{Thm: CompsExistRefined}
	Let $X$ be a curve with a unique non-Gorenstein singularity $p_0 \in X$. Then
	\begin{enumerate}
		\item the closure of $\nSm_d$ is a $g$-dimensional irreducible component of $\Jac^{d}$;\label{Item: FirstCompsExist}
		\item the closure of the inverse image $\AbelQ^{-1}(\nSm_{2 g -2 -d})$ is a $d$-dimensional irreducible component of $\Quot^{d}$ provided $d \ge 0$; \label{Item: SecondCompsExist}
		\item the closure of the inverse image $\AbelH^{-1}(\nSm_{-d})$ is a $n_{d}$-dimensional irreducible component of $\Hilb^{d}$ provided $d \ge d_0-g+1$. \label{Item: ThirdCompsExist}
	\end{enumerate}
	Furthermore, $\Jac^{d}$ (resp. $\Quot^{d}$, $\Hilb^{d})$ is $k$-smooth at a general closed point of $\nSm_{d}$ (resp. $\AbelQ^{-1}(\nSm_{2g-2-d})$, $\AbelH^{-1}(\nSm_{-d})$).
\end{theorem}

\begin{proof}  We give three separate arguments, one for each moduli space.

	\subsection*{The compactified Jacobian} The author claims that $\Jac^{d}$ is smooth of local dimension $g$ at every point $[I] \in \nSm_d$.  First, we establish this for $\Jac^{2g -2}$ at the point $[\dual]$ using deformation theory.  The groups $\ShExt^{q}(\dual, \dual)$ vanish for $q>0$ by \cite[Thm.~3.3.10]{bruns}.  We can conclude that the edge map $H^{p}(X, \calO_{X}) \to \Ext^{p}(\dual, \dual)$ associated to the local-to-global spectral sequence $H^{p}(\ShExt^{q}(\dual, \dual)) \Rightarrow \Ext^{p+q}(\dual, \dual)$ is an isomorphism for all $p$.  
	
	In particular, $\Ext^{2}(\dual, \dual)=0$ and $\Ext^{1}(\dual, \dual)$ is $g$-dimensional.  The vanishing of $\Ext^{2}(\dual, \dual)$ implies that $\Jac^{2g -2}$ is smooth at $[\dual]$, so the local dimension of $\Jac^{2g-2}$ at $[\dual]$ equals the tangent space dimension, $\dim \Ext^{1}(\dual, \dual) = g$.  This proves the claim for $I = \dual$, and we can conclude that the result holds for arbitrary $[I] \in \nSm_d$ by homogeneity. Thus the closure of $\nSm_{d}$ must be an irreducible component of $\Jac^{d}$ because the local dimension of $\overline{J}^{d}$ at any $[I] \in \nSm_{d}$ equals the dimension of $\nSm_{d}$.

\subsection*{The Quot scheme} When $d \ge 2g-1$, the Abel map $\AbelQ \colon \Quot^{d} \to \Jac^{2g-2-d}$ is smooth with fibers isomorphic to $\bbP^{d-g}$ \cite[Thm.~8.4(v)]{altman80}.  We can conclude from Item~(1) that the claim holds and additionally $\Quot^{d}$ is $k$-smooth at every closed point of $\AbelQ^{-1}(\nSm_{2g-2-d})$.  In fact, the proof of \cite[Thm.~8.4(v)]{altman80} shows that the obstruction group $\Ext^{1}(\ker(q), Q)$ vanishes for all $[q \colon \dual \twoheadrightarrow Q] \in \AbelQ^{-1}(\nSm_{2g-2-d})$.

We can deduce the case where $d$ is arbitrary from the case where $d$ is large.  Given a quotient map $q \colon \dual \twoheadrightarrow Q$ and  a collection of closed points $p_1, \dots, p_e \in X^{\text{sm}}$ disjoint from the support of $Q$,  we write
\begin{displaymath}
	q_{i} \colon \dual \twoheadrightarrow  \dual/\dual(-p_i) = k(p_{i})
\end{displaymath}
for the quotient map $\dual \twoheadrightarrow \dual/\dual(-p_i)$
and 
\begin{displaymath}
	q \times q_1 \times \dots \times q_{e} \colon \dual \twoheadrightarrow  Q \times k(p_1) \times \dots \times k(p_e)
\end{displaymath}

for the map into the product.  The locus $\AbelQ^{-1}(\nSm_{2g-2-d})$ is always non-empty, for it contains the points [$q_1 \times \dots \times q_d]$, $p_1, \dots, p_d \in X^{\text{sm}}$.  In fact, the closure $\overline{Y}_{d} \subset \Quot^{d}$ of the subset of all such points is irreducible and $d$-dimensional.  

$d$ is also equal to the dimension of the tangent space to $\Quot^{d}$  at any closed point of $\AbelQ^{-1}(\nSm_{2g-2-d})$.  We prove this as follows. Given $[q] \in \AbelQ^{-1}(\nSm_{2g-2-d})$,  fix $e$ large and $p_1, \dots, p_e$ general.  Then $[q' := q \times q_1 \times \dots \times q_e]$ lies in $\AbelQ^{-1}(\nSm_{-d-e})$, and the dimension of the tangent space at this point is $d+e$.  But this tangent space can be rewritten as

\begin{align*}
	\operatorname{T}_{[q']} \Quot^{d+e} 	=& \Hom(\ker(q'),Q \times k(p_1) \times \dots \times k(p_e))\\
						=& \Hom(\ker(q),Q) \oplus  \Hom(\ker(q_1), k(p_1) ) \oplus \cdots \\
						& \phantom{\oplus} \oplus \Hom(\ker(q_e), k(p_e))\\
						=& \operatorname{T}_{[q]} \Quot^{d} \oplus \Hom(\ker(q_1),k(p_1)) \oplus \cdots \\
						& \phantom{\oplus} \oplus \Hom(\ker(q_e),k(p_e)).
\end{align*}

Taking dimensions, we get $\dim \operatorname{T}_{[q]} \Quot^{d} = d$.  A similar computation shows that the $\Ext^{1}(\ker(q),Q)=0$, so $\Quot^{d}$ is $k$-smooth of local dimension $d$ at $[q]$.

We can conclude that $\overline{Y}_{d}$ is an irreducible component of $\Quot^{d}$.  To complete the proof, we must show that $\overline{Y}_{d}$ equals the closure of  $\AbelQ^{-1}(\nSm_{2g-2-d})$.  Thus suppose that $\overline{C}_{d}$ is an irreducible component of the closure.  Because the local dimension of $\Quot^{d}$ at any point of $\AbelQ^{-1}(\nSm_{2g-2-d})$ is $d$, the component $\overline{C}_{d}$ must have dimension $d$.  Now fix $e$ large.  If we define $\overline{C}_{d+e}$ to be the closure of the subset of points $[q \times q_1 \times \dots \times q_e]$ with $[q] \in \overline{C}_{d}$ and $p_1, \dots, p_e \in X^{\text{sm}}$ general, then $\overline{C}_{d+e}$ is contained in the closure of $\AbelQ^{-1}(\nSm_{2g-2-d-e})$.  Both subsets of $\Quot^{d+e}$ are $(d+e)$-dimensional, irreducible, and closed, hence the containment is an equality.  But $\overline{C}_{d+e}$ does not contain the general element of the form $[q_1 \times \dots \times q_{d+e}]$ (as the analogous statement holds for $\overline{C}_{d}$).  A contradiction! This completes the proof.

\subsection*{The Hilbert scheme} To begin, I claim that $\AbelH^{-1}(\nSm_{-d})$ is open in $\Hilb^{d}$.  It is, of course, enough to show that $\nSm_{-d} \subset \Jac^{-d}$ is open, and this can be established as follows.   Lemma~\ref{Lemma: LocClosed} states that $\nSm_{-d}$ is open in its closure $\overline{\nSm}_{-d}$.  Furthermore, the proof of Theorem~\ref{Thm: CompsExistRefined} shows that $\Jac^{-d}$ is smooth at every point of $\nSm_{-d}$, so the only  component of $\Jac^{-d}$ that meets $\nSm_{-d}$ is the closure $\overline{\nSm}_{-d}$.  We can conclude that $\nSm_{-d}$ is open in $\Jac^{-d}$ as its complement is the union of the irreducible components distinct from  $\overline{\nSm}_{-d}$ together with the closed subset $\overline{\nSm}_{-d} \setminus \nSm_{-d}$.  

Now consider  $\AbelH^{-1}(\nSm_{-d}^{\operatorname{o}}$.  This subset must also be open in $\Hilb^{d}$ as  $\nSm_{-d}^{\operatorname{o}}$ is open in $\nSm_{-d}$ by Lemma~\ref{Lemma: AbelForHilb}.  Furthermore, the same lemma states that  $\AbelH \colon \AbelH^{-1}(\nSm_{-d}^{\operatorname{o}}) \to \nSm_{-d}^{\operatorname{o}}$ is smooth with fibers isomorphic to $\bbP^{n_{d}-g}$.   We can immediately conclude that $\AbelH^{-1}(\nSm_{-d}^{\operatorname{o}})$ is $n_{d}$-dimensional, irreducible, and $k$-smooth.  The Hilbert scheme $\Hilb^{d}$ must also be $k$-smooth at every point of $\AbelH^{-1}(\nSm_{-d}^{\operatorname{o}})$ as $\AbelH^{-1}(\nSm_{-d}^{\operatorname{o}}) \subset \Jac^{-d}$ is open.  We can  conclude that  the closure of $\AbelH^{-1}(\nSm_{-d})$ is an irreducible component of $\Jac^{-d}$.  (The subset $\AbelH^{-1}(\nSm_{-d})$ is open in any irreducible component containing it, hence $\AbelH^{-1}(\nSm_{-d})$ is dense is any such component.)  This completes the proof.
\end{proof}

\section{Proof of Theorem~\ref{Thm: AllComps}}  \label{Sec: Module}
Here we prove Theorem~\ref{Thm: RefinedAllComps}, which is Theorem~\ref{Thm: AllComps} and the analogous statement for $\Quot^{d}$.  The theorem concerns non-planar curve singularities of finite representation type.  The classification of these singularities is recalled in Section~\ref{Sec: Tables}, where the singularities are listed in Table~\ref{Table: CurveSing}.  Table~\ref{Table: BigTable} of that section contains a list of the rank $1$, torsion-free modules over the ring of a singularity from Table~\ref{Table: CurveSing}.  We advise the reader to look at Section~\ref{Sec: Tables} before reading the proof of Theorem~\ref{Thm: RefinedAllComps}.

In proving Theorem~\ref{Thm: RefinedAllComps}, we need the following lemma, which we use to argue that it is enough to work with modules over $\widehat{\calO}_{X, p_0}$ rather
than sheaves over $X$.

\begin{lemma} \label{Lemma: LocToGlo}
	Let $X$ be a curve with a unique singularity $p_0 \in X$.  Suppose that $I$ is a rank $1$, torsion-free sheaf on $X$ and $\widehat{I}_{a} \subset  \widehat{\calO}_{X,p_0} \otimes k[[a]]$ is an ideal with  $k[[a]]$-flat quotient such that  there exists an isomorphism
	\begin{displaymath}
		\widehat{I}_{a} \otimes k[[a]]/(a) \cong I \otimes \widehat{\calO}_{X, p_0}.
	\end{displaymath}
	
	Then there exists a $k[[a]]$-flat family of rank $1$, torsion-free sheaves $I_{a}$ on $X \otimes k[[a]]$ with the property that there exists an isomorphism
	\begin{displaymath}
		\widehat{I}_{a} \cong I_{a} \otimes ( \widehat{\calO}_{X, p_0} \otimes k[[a]] ).
	\end{displaymath}
\end{lemma}
\begin{proof}
Define $J_{a} \subset \calO_{X} \otimes k[[a]]$ to be the kernel of the composition 
\begin{displaymath}
	\calO_{X} \otimes k[[a]]  \to \widehat{\calO}_{X, p_0} \otimes k[[a]] \to  \widehat{\calO}_{X, p_0} \otimes k[[a]] /\widehat{I}_{a}.
\end{displaymath}

The quotient $\calO_{X} \otimes k[[a]]/J_{a}$ is canonically isomorphic to $\widehat{\calO}_{X, p_0} \otimes k[[a]]/\widehat{I}_{a}$, and hence is $k[[a]]$-flat. In particular, $J_{a}$ is itself $k[[a]]$-flat.  Furthermore, the fibers of $J_{a}$ are  rank $1$, torsion-free sheaves.  Indeed, the generic fiber of $J_{a}$ is rank $1$ and torsion-free because it is a subsheaf of $\calO_{X} \otimes \operatorname{Frac} k[[a]]$ that is nonzero (as the quotient  is supported at $p_0$).  Similarly, because the quotient $\calO_{X} \otimes k[[a]]/J_{a}$ is $k[[a]]$-flat, the reduction $J := J_{a} \otimes k[[a]]/(a) \to \calO_{X}$ of the inclusion map is injective, and so the special fiber $J$ is a nonzero subsheaf of $\calO_{X}$, hence rank $1$ and torsion-free.  We have now  shown that  the sheaf $J_{a}$ has all of the desired properties except that the special fiber $J := J_{a} \otimes k[[a]]/(a)$ may not be isomorphic to $I$.  

We proceed to modify $J_{a}$ so that $J$ is isomorphic to $I$.  While $J$ and $I$ may not be isomorphic, these two sheaves do become isomorphic after  passing to $\widehat{\calO}_{X, p_0}$.  Thus the completed stalk of $\ShHom(I,J)$ at $p_0$ is free of rank $1$ as the formation of  $\ShHom(I,J)$ commutes with completion.  Consequently, there exists an open neighborhood $U \subset X$ of $p_0$ and an isomorphism $\phi_{1} \colon J|_{U} \cong I|_{U}$. (Pick $\phi_1$ to map to a generator of the completion of $\ShHom(I,J)$.)  Away from $p_0$, the sheaves $J$ and $I$ are locally isomorphic because both sheaves restrict to line bundles on $X \setminus \{ p_0 \}$.  Now the complement $X \setminus U$ consists of a finite number of points, so we can find an open subset  $V \subset X$ that contains $X \setminus U$ and has the property that there exists an isomorphism $\phi_{2} \colon J|_{V} \cong I|_{V}$.  On the overlap $U \cap V$, the automorphism $\phi_{2}^{-1} \circ \phi_{1} \colon J|_{U \cap V} \cong J|_{U \cap V}$ is an automorphism of a line bundle and so is defined by multiplication with a fixed function $f \in H^{0}(U \cap V, \calO_{X}^{\ast})$.  Define $L$ to be the line bundle obtained by glueing $\calO_{V}$ to $\calO_{U}$ over $U \cap V$ by the automorphism defined by $f$.  Then an isomorphism 
$$
	J \otimes L \cong I
$$
is defined by 
\begin{gather*}
	s \otimes 1 \mapsto \phi_{1}(s) \text{ on $U$,} \\
	s \otimes 1 \mapsto \phi_{2}(s) \text{ on $V$}.
\end{gather*}
We can conclude that the tensor product of $J_{a}$ with the constant family of line bundles with fiber $L$ satisfies all of the desired properties.  This completes the proof.   (We thank the anonymous referee for suggesting this argument.)
\end{proof}

We now prove the main theorem.
\begin{theorem} \label{Thm: RefinedAllComps}
	Let $X$ be a curve with a unique singularity $p_0 \in X$ that is of finite representation type and non-Gorenstein.  Then
	\begin{enumerate}
		\item $\Jac^{d}$ has exactly two irreducible components: the smoothable component and the closure of $\nSm_{d}$;
		\item $\Quot^{d}$ has exactly two irreducible components provided $d \ge 2g-1$: the smoothable component and closure of  $\AbelQ^{-1}(\nSm_{2g-2-d})$.
	\end{enumerate}
\end{theorem}
\begin{proof}
First, we reduce to the problem of deforming the modules in Table~\ref{Table: BigTable}, and then we deform those modules on a case-by-case basis.  To make the reduction, consider the following hypothesis:
\begin{hypothesis} \label{Hypothesis}
	Let $R$ be a ring that is a finite product of rings of curve singularities.  If $M$ is a rank $1$, torsion-free $R$-module, then we say that $(R,M)$ satisfies Hypothesis~\ref{Hypothesis} if there exists an ideal $M_{a} \subset R \otimes k[[a]]$ with $k[[a]]$-flat cokernel such that there exists an isomorphism
	\begin{displaymath}
		M_{a} \otimes k[[a]]/(a) \cong M
	\end{displaymath}
	over the special fiber and an isomorphism over the completed generic fiber that is either of the form
	$$
		M_{a} \widehat{\otimes} \operatorname{Frac} k[[a]] \cong R \widehat{\otimes} \operatorname{Frac} k[[a]]
	$$
	or of the form 
	$$
		M_{a} \widehat{\otimes} \operatorname{Frac} k[[a]] \cong \dual \widehat{\otimes} \operatorname{Frac} k[[a]].
	$$

\end{hypothesis}
Note that we require the existence of an isomorphism over the completed tensor product $R \widehat{\otimes} \operatorname{Frac} k[[a]]$, not the uncompleted tensor product  $R \otimes \operatorname{Frac} k[[a]]$.  The ring $R \otimes \operatorname{Frac} k[[a]]$ may fail to be complete, and it is the relevant completed tensor product that is isomorphic to the completed local ring of $X \otimes \operatorname{Frac} k[[a]]$ at $p_0$ for $X$ as in the statement of the theorem.  We also point out that $R$ is a finite product of rings of singularities, not the  ring of a singularity.  We allow for finite products because a finite extension of the ring of a singularity may be a product of rings of singularities.  (E.g.~$k[[t]] \times k[[t]]$ is a finite extension of the ring of the $A_1$-singularity.)

The theorem quickly follows if we assume Hypothesis~\ref{Hypothesis} holds when $R=\widehat{\calO}_{X, p_0}$ and $M$ is arbitrary.  Indeed, let us prove this first for the compactified Jacobian and then for the Quot scheme.

\subsection*{The compactified Jacobian} For Item~(1), we need to show that $J^{d}_{X} \cup \nSm_d$ is dense in $\Jac^{d}$.  It is enough to prove that the closure of $J^{d}_{X} \cup \nSm_d$ contains every closed point of $\Jac^{d}$, so let  $[I] \in \Jac^{d}$ be a given closed point.  Apply Hypothesis~\ref{Hypothesis} to $M := I \otimes \widehat{\calO}_{X,p_0}$.  If $\widehat{I}_{a} := M_{a}$ is as in the conclusion of Hypothesis~\ref{Hypothesis}, then Lemma~\ref{Lemma: LocToGlo} asserts that there is a flat deformation $I_{a}$ of $I$ with the property that 
\begin{displaymath}
	\widehat{I}_{a}  \cong I_{a} \otimes (\widehat{\calO}_{X, p_0} \otimes k[[a]]).
\end{displaymath}
The family $I_{a}$ corresponds to a morphism $\Spec(k[[a]]) \to \Jac^{d}$ that sends the special point to $[I]$ and the generic point to an element of $J^{d}_{X} \cup \nSm_{d}$, proving Item~(1) (under the assumption that Hypothesis~\ref{Hypothesis} holds).

\subsection*{The Quot scheme} Given Item~(1), Item~(2) follows immediately as $\AbelQ$ is a $\bbP^{d-g}$-bundle for $d \ge 2 g -1$.  

\subsection*{Hypothesis~\ref{Hypothesis} holds}  We now prove that Hypothesis~\ref{Hypothesis} holds when $R$ is of finite representation type. Our strategy is as follows. We begin by making some preliminary reductions.  These reductions will let us set up an inductive argument which reduces the claim that Hypothesis~\ref{Hypothesis} holds to the claim  that fifteen specific modules deform.  We complete the proof by deforming these modules by hand.

\subsubsection*{Reduction One}  In Hypothesis~\ref{Hypothesis}, we may assume that the rank $1$, torsion-free module $M$ is an ideal $M \subset R$. Indeed, we can construct an embedding as follows.  The natural map $M \to M \otimes \operatorname{Frac} R$ is injective (as $M$ is torsion-free) and $M \otimes \operatorname{Frac} R$ is isomorphic to $\operatorname{Frac} R$ (as $M$ is rank $1$).  If we fix an isomorphism $M \otimes \operatorname{Frac} R \cong \operatorname{Frac} R$, then the composition 
\begin{displaymath}
	M \stackrel{i}{\longrightarrow} M \otimes \operatorname{Frac} R \cong \operatorname{Frac} R
\end{displaymath}	
is an injection.  The image may not lie in $R$, but if we fix a nonzero divisor $t \in R$, then the image of $t^{b} \cdot i$ will lie in $R$ once $b$ is sufficiently large.  

\subsubsection*{Reduction Two} Hypothesis~\ref{Hypothesis} holds when $R$ is a finite product of rings of planar singularities.  It is enough to consider the case $R=\widehat{\calO}_{X, p_0}$ for $X$ a locally planar curve.  Let $M$ be given. We have just shown that we can realize $M$ as an ideal $M \subset \widehat{\calO}_{X, p_0}$.  Define $Z \subset X$ to be the closed subscheme of $X$ that corresponds to the quotient map $\calO_{X} \to \widehat{\calO}_{X,p_0} \to \widehat{\calO}_{X, p_0}/M$.  The main result of  \cite{altman77} implies that the point  $[Z]$ of the Hilbert scheme $\Hilb^{d}$ lies in the closure of the locus of Cartier divisors.  By \cite[Prop.~7.1.4]{EGA2}, this containment is witnessed by a morphism $S \to \Hilb^{d}$ out of the spectrum of a valuation ring that maps the special point to $[Z]$ and the generic point to a point in the locus of Cartier divisors.  Furthermore, $S$ can be chosen so that it is the spectrum of a complete discrete valuation ring with residue field $k$. Hence $S$ is isomorphic to $\Spec(k[[a]])$, and we obtain a suitable ideal $M_{a}$ by pulling back the universal family to $X \times S$ and then further restricting to $R \otimes k[[a]]$.

\subsubsection*{Reduction Three} We now establish a result that  will let us set up an induction.  Let $\calO = \widehat{\calO}_{X, p_0}$ be the completed local ring of the curve $X$ at a closed point $p_0$.  Suppose $\calO \subset \calO'$ is a  finite extension contained in $\operatorname{Frac} \calO$ and assume Hypothesis~\ref{Hypothesis} is satisfied in the following two cases: 
\begin{itemize}
	\item $R$ equals $\calO$ and $M$ equals the ring $\calO'$ or its dualizing module $\dual'$;
	\item $R$ equals $\calO'$ and $M$ equals an arbitrary $\calO'$-module $M'$.
\end{itemize}
Then Hypothesis~\ref{Hypothesis} holds when $R=\calO$ and $M=M'$ is a $\calO'$-module considered as a $\calO$-module.
To prove this, say $\widehat{I}' := M'$ is a given such module.

We can construct a finite birational morphism $f \colon X' \to X$ such that $\calO'$ is the completion of the localization of $X'$ at a finite set of closed points.  Indeed, $\calO'$ is a rank $1$, torsion-free module over $\calO$, so it is isomorphic to some ideal $\widehat{J} \subset \calO$.  If we form the kernel $J$ of the composition $\calO_{X} \to \calO \to \calO/\widehat{J}$, then we can take $X' = \Spec(\ShEnd(J))$.  The ring $\calO'$ is then the completion of the localization of $X'$ at $f^{-1}(p_0)$.  Similarly, we can assume $\widehat{I'}$ is isomorphic to $I' \otimes \calO'$ for some rank $1$, torsion-free sheaf $I'$ on $X'$.  (Construct $I'$ as an ideal.)

Now Hypothesis~\ref{Hypothesis} is satisfied when $R=\calO'$, so we can conclude that $J^{d'}_{X'} \cup \nSm'_{d'}$ is dense in $\overline{J}_{X'}^{d'}$.  By \cite[Prop.~3]{altman90}, this compactified Jacobian $\overline{J}_{X'}^{d'}$ embeds in $\overline{J}_{X}^{d}$ by the rule $I' \mapsto f_{*}(I')$ (for suitable $d$).   Now Hypothesis~\ref{Hypothesis} is also satisfied when $R=\calO'$ and $M=\calO'$ or $\dual'$, so we can conclude that the closure of $J^{d}_{X} \cup \nSm_d$ contains the image of $J^{d'}_{X'} \cup \nSm'_{d'}$, and hence the entire image of $\overline{J}_{X'}^{d'}$. 

We now construct a suitable deformation of $\widehat{I'}$ as follows.  We have just shown that $[f_* I']$ lies in the closure of $J^{d}_{X} \cup \nSm_d$, and this containment is witnessed by a morphism out of $S=\Spec(k[[a]])$ \cite[Prop.~7.1.4]{EGA2}.  Pulling back the universal family to $\widehat{\calO}_{X, p_0} \otimes k[[a]]$ produces a suitable module $M_{a}=\widehat{I}_{a}$ except that $\widehat{I}_{a}$ is not obviously an ideal.

We can, however, arrange that $\widehat{I}_{a}$ is an ideal as follows.  We can make the degree $d$ of $I$ as large as we wish, and if we make the degree  large enough, then the Abel map $\AbelQ$ out of $\Quot^{d}$ is smooth.  In particular, we can lift  $S \to \Jac^{d}$ to a morphism $S \to \Quot^{d}$, and we can thus assume that $\widehat{I}_{a} \subset \dual \otimes k[[a]]$ is a submodule with $k[[a]]$-flat cokernel.  If we fix an injection $\dual \hookrightarrow \calO$, then the composition 
\begin{displaymath}
	\widehat{I}_{a} \subset \dual \otimes k[[a]] \hookrightarrow \calO \otimes k[[a]]
\end{displaymath}
realizes $\widehat{I}_{a}$ as a suitable ideal.

\subsubsection*{Hypothesis~\ref{Hypothesis} holds}
We now prove that Hypothesis~\ref{Hypothesis} holds when $R=\calO$ is the ring of  a  non-Gorenstein singularity that is of finite representation type.  Our argument makes use of the classification of rank $1$, torsion-free modules over such a ring, and it is recommended that the reader look at Section~\ref{Sec: Tables} before proceeding.

Table~\ref{Table: BigTable} of Section~\ref{Sec: Tables}  lists the modules over the ring of a non-Gorenstein curve singularity that is of finite representation type.  The rank $1$, torsion-free $\calO$-modules are all listed together.  The portion of the table containing the modules over a fixed $\mathcal{O}$ is further subdivided by horizontal bars.  (E.g.~the modules $R+R\cdot t^{2}$ and $R+R\cdot t$ over the $E_{8}(1)$-singularity are in the same subdivision.)  The modules within a given subdivision are arranged so that the endomorphism ring of a module $M$ is contained in the endomorphism ring of the module directly below it.  (E.g.~for the $E_{8}(1)$-singularity, the endomorphism ring of $R+R \cdot t^{2}$ is contained in the endomorphism ring of $R+R \cdot t$.)

We now use this table to prove that Hypothesis~\ref{Hypothesis} holds when $R=\calO$ is the ring of a non-Gorenstein singularity that is of finite representation.  We induct on the delta invariant $\delta(\calO) = \dim_{k} \widetilde{\calO}/\calO$.  Thus let $\calO$ be non-Gorenstein and of finite representation type and assume Hypothesis~\ref{Hypothesis} holds whenever $R$ is a ring that also satisfies these conditions but has delta invariant strictly smaller than $\delta(\calO)$.

By the induction, it is enough to prove that the hypothesis is satisfied when $M=\widehat{I}$ equals a module $\widehat{I}$ that is the topmost element of  a subdivision of Table~\ref{Table: BigTable} that is not equal to the ring $\calO$ or its dualizing module $\dual$.  (E.g.~when $\calO$ is the $E_{8}(1)$-singularity, it is enough to show the hypothesis is satisfied by $M=R+R \cdot t$ and $M=R+R \cdot t^{4}$.) Indeed, suppose $\widehat{I}$ is any rank $1$, torsion-free $\calO$-module.  If $\widehat{I}$ is equal to $\calO$, $\dual$, or one of the modules that we are assuming satisfies Hypothesis~\ref{Hypothesis}, then there is nothing to prove.  Otherwise, let $\widehat{I}_0$ be the topmost module in the subdivision containing $\widehat{I}$ that does not equal $\calO$ or $\dual$.  Set $\calO' := \End(\widehat{I})$ and $\calO'_0 := \End(\widehat{I}_0)$.  Inspecting Table~\ref{Table: BigTable}, we see that $\calO \subsetneq \calO'_0 \subsetneq \calO'$.  
 
We now use the inductive hypothesis. The ring $\calO'_0$ is a ring of finite representation type and satisfies $\delta(\calO'_0) < \delta(\calO)$.  By applying either  the inductive hypothesis (when $\calO'_0$ is non-Gorenstein) or Reduction Two (when $\calO'_0$ is planar), we can conclude that Hypothesis~\ref{Hypothesis} is satisfied when $R=\calO'_0$ and $M=\widehat{I}$ (considered as a $\calO'_0$-module).  Furthermore, an inspection of Table~\ref{Table: BigTable} shows that Hypothesis~\ref{Hypothesis} is satisfied by assumption when $R=\mathcal{O}$ and $M$ equals either the ring $\calO'_0$ and its dualizing module $\dual'_0$ (i.e.~Such an $M$ is topmost within its subdivision).  We can therefore conclude that Hypothesis~\ref{Hypothesis} holds when $R=\calO$ and $M=\widehat{I}$ by Reduction Three.  This proves the assertion.  We now complete the proof by showing that Hypothesis~\ref{Hypothesis} holds when $M = \widehat{I}$ is a topmost element.

There are fifteen such modules, and we deform these modules one-by-one.  In every case, the general technique is the same.  Given one of the fifteen modules $\widehat{I}$, we begin by constructing a surjective $\calO \otimes k[[a]]$-linear map 
\begin{equation} \label{Eqn: ModelSurj}
	\phi_{a} \colon \widehat{P} \otimes k[[a]] \to Q \otimes k[[a]].
\end{equation}
Here $\widehat{P}$ is a rank $1$, torsion-free module and $Q$ is a module of finite length.

In the constructions below, it is perhaps not always clear what the $\calO \otimes k[[a]]$-module structure  on $Q \otimes k[[a]]$ is.  Most often $Q=k$, and the $\calO$-module structure is defined by making $f \in \calO$ act as $f \cdot v = f(0) v$ for $v \in k$.  The $\calO \otimes k[[a]]$-module structure on $Q \otimes k[[a]]]$ is then defined by extending scalars.  In a few constructions, however, $Q=k[\epsilon]/(\epsilon^{2})$ or $k[\epsilon_1, \epsilon_2](\epsilon_1, \epsilon_2)^2$.  These $k$-algebras do not have distinguished  $\calO$-module structure, but the map $\phi_{a}$ that we construct will be a $k[[a]]$-algebra map.  In this case, we endow $Q \otimes k[[a]]$ with the induced $\calO \otimes k[[a]]$-module structure.

The kernel $\widehat{I}_{a} := \ker(\phi_{a})$ certainly has $k[[a]]$-flat cokernel --- the cokernel is $Q \otimes k[[a]]$.  The conclusion of Hypothesis~\ref{Hypothesis} requires that $\widehat{I}_{a}$ is an ideal, and we can arrange this by fixing an injection $\widehat{P} \hookrightarrow \calO$ and then forming the composition $I_{a} \hookrightarrow \widehat{P} \otimes k[[a]] \hookrightarrow \calO \otimes k[[a]]$.

To verify that $\widehat{I}_{a}$ satisfies the conclusion of Hypothesis~\ref{Hypothesis}, we also need to  exhibit an isomorphism
$$
	\widehat{I}_{a} \otimes k[[a]]/(a)	\cong \widehat{I}
$$
and either an isomorphism
$$
	\widehat{I}_{a} \widehat{\otimes} \operatorname{Frac} k[[a]] \cong \calO \widehat{\otimes} \operatorname{Frac} k[[a]]
$$
or an isomorphism 
$$
		\widehat{I}_{a} \widehat{\otimes} \operatorname{Frac} k[[a]] \cong \dual \widehat{\otimes} \operatorname{Frac} k[[a]]
$$
In every construction, both  $\widehat{I}$ and $\widehat{P}$ are submodules of $\operatorname{Frac} \calO$, and the isomorphism $\widehat{I}_{a} \otimes k[[a]]/(a) \cong \widehat{I}$ is constructed as the restriction of the map $\operatorname{Frac} \calO \to \widetilde{\calO}$ given by multiplication with a fixed nonzero divisor $f \in \widetilde{\calO}$.  Such a map is always injective, so we just need to check that the map is well-defined and surjective.  The second isomorphism is constructed in a similar manner.

We construct these isomorphisms and  the surjection $\phi_{a}$ using power series methods.  The normalization $\widetilde{\calO}$ is isomorphic to a self-product of the power series ring $k[[t]]$.  The tensor product $k[[t]] \otimes k[[a]]$ is not isomorphic to $k[[t,a]]$, but the natural map $k[[t]] \otimes k[[a]] \to k[[t,a]]=k[[a]][[t]]$ is injective.  Thus we can think of an element of the tensor product as a power series in $t$ and $a$, or alternatively as a power series in $t$ with coefficients in $k[[a]]$.  Given $f \in k[[t]] \otimes_{k} k[[a]]$, we write $f_n \in k[[a]]$ for the coefficient of $t^n$ in the image of $f$ in $k[[a]][[t]]$.  We will also abuse notation and write $a f$ in place of $f \otimes a$.  

We now proceed to show that the fifteen modules from Table~\ref{Table: BigTable} satisfy Hypothesis~\ref{Hypothesis}.  We will construct all of the relevant maps, but we do not always verify that maps have the desired properties (e.g.~that $\phi_{a}$ is a surjective homomorphism).  We verify these details for the first module only and  leave remaining cases to the interested reader.  

To make notation consistent with that of Table~\ref{Table: BigTable} (and \cite{greuel85}), we will write ``$R$" in place of ``$\calO$" and ``$M$" in place of ``$\widehat{I}$" for the remainder of this section.

\subsubsection*{The $A_{n} \vee L$-singularity, $n$ even} We need to deform the three modules $R + R \cdot (t^{n-1},0)$,
$R \cdot (1,0) + R \cdot (t^{n-3}, 1)$, and $R + R \cdot (1,0)$. 

We begin by deforming  $R + R \cdot (t^{n-1},0)$ to the dualizing module $\dual = R \cdot (1,0) + R \cdot (t^{n-1},1)$.  Define
\begin{gather*}
	\phi_{a} \colon  R + R \cdot (1,0)+R \cdot (t^{n-1}, 0) \otimes k[[a]] \to k[[a]], \\
	(f,g)  \mapsto f_{0} - g_{0} - a \left( f_{n-1}-g_{0} \right)
\end{gather*}
To see that $\phi_{a}$ is a $R \otimes k[[a]]$-module homomorphism map, observe that a typical element of the source can be written as a pair $(f,g)$ of power series satisfying $f_1=f_3=\dots=f_{n-3}=0$, and such a power series lies in $R$ when we additionally have $f_{n-1}=0$ and $f_{0}=g_{0}$.  Given $(p,q) \in R$ and $(f,g) \in R + R \cdot (1,0) + R \cdot (t^{n-1},0)$, we compute:
\begin{align*}
	\phi_{a}( (p,q) \cdot (f,g) ) =& p_{0} f_{0} - q_{0} g_{0}- a \left( p_0 f_{n-1}-q_0 g_0 \right)\\
				=& p_{0} f_{0} - p_{0} g_{0}- a \left( p_0 f_{n-1}-p_0 g_0 \right)\\
				=& p_0 \phi_{a}(f,g) \\
				=& (p,q) \cdot \phi_{a}(f,g).
\end{align*}
This shows that $\phi_a$ is linear, and the map is visibly surjective.

The fiber $M_{a} \otimes k[[a]]/(a)$ is equal to the submodule $R + R \cdot (t^{n-1},0) \subset \widetilde{R}$, so to show that $M_{a}$ satisfies
the conclusion of Hypothesis~\ref{Hypothesis}, we just need to verify that the completed generic fiber is isomorphic to $\dual \widehat{\otimes} \operatorname{Frac} k[[a]]$.  An 
isomorphism between these modules is given by
\begin{gather}
	 M_{a} \widehat{\otimes} \operatorname{Frac} k[[a]] \to 	\dual \widehat{\otimes} \operatorname{Frac} k[[a]], \label{Eqn: CpltGenericIso} \\
	(f,g) \mapsto  (a-t^{n-1}, a-1) \cdot (f,g). \notag
\end{gather}
This map is well-defined because 
\begin{align*}	
	(a-t^{n-1}, a-1) \cdot (f,g) =& \left( (a-t^{n-1})f- (a-1) g_0 t^{n-1}, a f_0 \right) \cdot (1,0)+\\
		&  \phantom{( (a-t^{n-1})f- } \left( (a-1) g_0, (a-1) g \right) \cdot (t^{n-1}, 1),
\end{align*}
and an inspection of the relevant power series shows  that the coefficients appearing in the right-hand side of the above equation lie in
$R \widehat{\otimes} \operatorname{Frac} k[[a]]$ provided $(f,g) \in M_{a} \widehat{\otimes} \operatorname{Frac} k[[a]]$.

The injectivity of Eq.~\eqref{Eqn: CpltGenericIso} is automatic, and surjectivity follows from the identities
\begin{gather*}
	(a-t^{n-1},a-1) \cdot (a^{-1} + a^{-2} t^{n-1} + \dots,0)=(1,0),\\
	(a-t^{n-1},a-1) \cdot (a^{-1} t^{n-1}+ a^{-2} t^{2(n-1)} + \dots, (a-1)^{-1})=(t^{n-1},1).
\end{gather*}
This shows that the conclusion of Hypothesis~\ref{Hypothesis} holds when $M=R + R \cdot (t^{n-1},0) \subset \widetilde{R}$.

We now deform the module $R \cdot (1,0) + R \cdot (t^{n-3},1)$ to $R$.  Define 
\begin{gather*}
	\phi_{a} \colon R \otimes k[[a]] \to k[\epsilon_1, \epsilon_2]/(\epsilon_1, \epsilon_2)^{2} \otimes k[[a]], \\
	(f,g) \mapsto f_{0} + \epsilon_1 \left(f_{2} +  a f_{n+1} \right) + \epsilon_2 \left( g_{1} + (a-1) f_{n+1} \right).
\end{gather*} 
A computation shows that this is a surjective $k[[a]]$-algebra map.

The kernel $M_{a} := \ker \phi_{a}$ is a deformation of $R \cdot (1,0) + R \cdot (t^{n-3},1)$ because the homomorphism 
\begin{gather*}
	R \cdot (1,0) + R \cdot (t^{n-3},1) \to M_{a} \otimes k[[a]]/(a), \\
	(f,g) \mapsto (t^{4}, t) \cdot (f,g)
\end{gather*}
is an isomorphism.  The module is a deformation to $R$ because the map
\begin{gather*}
	R \widehat{\otimes} \operatorname{Frac} k[[a]] \to M_{a} \widehat{\otimes} \operatorname{Frac} k[[a]], \\
	(f,g) \mapsto (t^{n+1}-a t^{2}, (1-a) t) \cdot (f,g)
\end{gather*}
is an isomorphism.  This shows that  Hypothesis~\ref{Hypothesis} is satisfied when $M=R \cdot (1,0) + R \cdot (t^{n-3},1)$.

Finally, we show that  $R + R \cdot (1,0)$ deforms to $R$. Define
\begin{gather*}
	\phi_{a} \colon R + R \cdot (1,0) \otimes k[[a]] \to k[[a]], \\
	(f,g) \mapsto g_{0} - a f_{0} - a g_{0}.
\end{gather*}
As with the previous maps, a computation shows that $\phi_{a}$ is a $R \otimes k[[a]]$-linear surjection.  Set $M_{a} = \ker \phi_{a}$.
Because the maps 
\begin{gather*}
	R + R \cdot (1,0) \to M_{a} \otimes k[[a]]/(a), \\
	(f,g) \mapsto (1,t) \cdot (f,g)
\end{gather*}
and
\begin{gather*}
	M_{a} \widehat{\otimes} \operatorname{Frac} k[[a]] \to R \widehat{\otimes} \operatorname{Frac} k[[a]], \\
	(f,g) \mapsto (a, (1-a)) \cdot (f,g).
\end{gather*}
are isomorphisms, $M_{a}$ is a deformation of $R+ R \cdot (1,0)$ to $R$.  This proves that  Hypothesis~\ref{Hypothesis} is satisfied when $M=R + R \cdot (1,0)$.

\subsubsection*{The $A_{n} \vee L$-singularity, $n$ odd}  There are five modules that we need to deform:
$R + R \cdot (t^{(n-1)/2},0,0)$, $R \cdot (1,1,0) + R \cdot (t^{(n-3)/2},0,1)$,
$R + R \cdot (1,0,0)$, $R + R \cdot (0,1,0)$, and $R +R \cdot (0,0,1)$.  We begin with the first module.

We deform $R + R \cdot (t^{(n-1)/2},0,0)$ to $\dual=R \cdot (1,1,0) + R \cdot (t^{(n-1)/2},0,1)$.  Define 
\begin{gather*}
	\phi_{a} \colon R \otimes k[[a]] \to k[[\epsilon]]/(\epsilon^{2}) \otimes k[[a]], \\
	(f,g,h) \mapsto f_0 + \epsilon \left( f_1+g_1-2h_1 - a f_{(n+1)/2} + a g_{(n+1)/2}+ a h_{1} \right).
\end{gather*}
This is a surjective homomorphism.  Set $M_{a} := \ker \phi_{a}$.  The homomorphisms 
\begin{gather*}
	R + R \cdot (t^{(n-1)/2},0,0) \to M_{a} \otimes k[[a]]/(a), \\
	(f,g,h) \mapsto (t,t,t) \cdot (f,g,h),
\end{gather*}
and
\begin{gather*}
	M_{a} \widehat{\otimes} \operatorname{Frac} k[[a]] \to \dual \widehat{\otimes} \operatorname{Frac} k[[a]], \\
	(f,g,h) \mapsto (a t^{-1}-t^{(n-3)/2}, a t^{-1}+t^{(n-3)/2},a t^{-1} -2 t^{-1}) \cdot (f,g,h)
\end{gather*}
are isomorphisms, so the module $M_{a}$ satisfies the conclusion of Hypothesis~\ref{Hypothesis}.

A similar construction deforms the module $R \cdot (1,1,0) + R \cdot (t^{(n-3)/2},0,1)$ to $R$.  Form the kernel $M_{a}$ of the map 
\begin{gather*}
	\phi_{a} \colon R \otimes k[[a]] \to k[\epsilon_1, \epsilon_2]/(\epsilon_1, \epsilon_2)^2 \otimes k[[a]], \\
	(f,g,h) \mapsto f_{0} + \epsilon_{1} \left( (1+a) (f_{(n+1)/2} - g_{(n+1)/2}) - h_1\right) \\
		\phantom{(f,g,h) \mapsto a_{0} + \epsilon_{1}(1+a) (f_{(n+1)/2}}		{}+ \epsilon_2 \left(  (1+a) (f_1+g_1 ) -  a h_1 \right).
\end{gather*}
As before, this map is a surjective homomorphism.  The rule
\begin{gather*}
	R \cdot (1,1,0) + R \cdot (t^{(n-3)/2},0,1) \to M_{a} \otimes k[[a]]/(a), \\
	(f,g,h) \mapsto (t^{2}, t^{2}, t) \cdot (f,g,h)		
\end{gather*}
defines an isomorphism, which shows that $M_{a}$ is a deformation of $R \cdot (1,1,0) + R \cdot (t^{(n-3)/2},0,1)$.  This deformation is a deformation to $R$ because the rule 
\begin{gather*}
	R \widehat{\otimes} \operatorname{Frac} k[[a]] \to 	M_{a} \widehat{\otimes} \operatorname{Frac} k[[a]], \\
	(f,g,h) \mapsto (a +t^{(n-1)/2}, a-t^{(n-1)/2},1+a) \cdot (t,t,2t) \cdot (f,g,h).
\end{gather*}
defines an isomorphism.  We can conclude that Hypothesis~\ref{Hypothesis} is satisfied when $M=R \cdot (1,1,0) + R \cdot (t^{(n-3)/2},0,1)$.

Next we deform $R + R \cdot (1,0,0)$ to $R$.  One such deformation is the kernel $M_{a}$ of the
homomorphism 
\begin{gather*}
	\phi_{a} \colon R + R \cdot (1,0,0) \otimes k[[a]] \to k[t]/(t^{(n+1)/2}) \otimes k[[a]], \\
	(f,g,h) \mapsto (1+a) f(t) - a g(t) \text{  (mod $t^{(n+1)/2})$,}
\end{gather*}
which is a surjective homomorphism.  The maps 
\begin{gather*}
	R + R \cdot (1,0,0) \to M_{a} \otimes k[[a]]/(a), \\
	(f,g,h) \mapsto (t^{(n-1)/2}, 1, 1) \cdot (f,g,h)
\end{gather*}
and
\begin{gather*}
	M_{a} \widehat{\otimes} \operatorname{Frac} k[[a]] \to R \widehat{\otimes} \operatorname{Frac} k[[a]] , \\
	(f,g,h) \mapsto (a^{-1}-1,1,1) \cdot (f,g,h)
\end{gather*}
are isomorphisms, and this shows that  $M=R + R \cdot (1,0,0)$ satisfies Hypothesis~\ref{Hypothesis} .

To complete this case, we need to deform the modules $R + R \cdot (0,1,0)$  and $R + R \cdot (0,0,1)$.
If we modify the  construction of the deformation of $R + R \cdot (1,0,0)$  by swapping the roles of $f$ and
$g$, then we obtain a suitable deformation $R + R \cdot (0,1,0)$ to $R$.  If we instead swap $f$ and $h$ and change
the target of $\phi_{a}$ to $k[[a]]$, then we obtain a deformation of  the module  $R + R \cdot (0,0,1)$ to $R$.

\subsubsection*{The $E_6(1)$-singularity}  It is enough to show that $R+R \cdot t^2$ deforms to the dualizing module $\dual = R + R \cdot t$.  Define $M_{a}$ to be the kernel of the map
\begin{gather*}
	\phi_{a} \colon \widetilde{R} \otimes k[[a]] \to k[[a]], \\
	f \mapsto f_1-a f_2.
\end{gather*}
This map is a  surjective  homomorphism.  The fiber $M_{a} \otimes k[[a]]/(a)$ of the kernel $M_{a} := \ker(\phi_{a})$  is equal to $R+R \cdot t^2$, and the map 
\begin{gather*}
	M_{a} \widehat{\otimes} \operatorname{Frac} k[[a]] \to R+R\cdot t \widehat{\otimes} \operatorname{Frac} k[[a]], \\
	f \mapsto (1-a^{-1}t) \cdot f
\end{gather*}
is an isomorphism.  This proves that $M=R+R\cdot t^{2}$ satisfies Hypothesis~\ref{Hypothesis}.

\subsubsection*{The $E_{7}(1)$-singularity}  The endomorphism ring of $R + R \cdot (1,0)$ is contained in the endomorphism ring of $R + R \cdot (t^2,0)$, so it is enough to deform $R + R \cdot (t^{2}, 0)$ and $R \cdot (1,0) + R \cdot (t,1)$.  

The module $R + R \cdot (t^{2}, 0)$ deforms to the dualizing module $\dual=R+R\cdot (t,0)$.  Indeed, the map
\begin{gather*}
	\phi_{a} \colon R + R \cdot (t,0) + R \cdot (t^{2}, 0) \otimes k[[a]] \to k[[a]], \\
	(f,g)  \mapsto f_1 - a  \left(f_2 - g_1\right).
\end{gather*}
is a linear surjection, and the kernel $M_{a} = \ker \phi_{a}$ is a suitable deformation. To see this, observe that the fiber $M_{a} \otimes k[[a]]/(a)$ is equal to $R+R\cdot (t^{2},0)$, and the map 
\begin{gather*}
	M_{a} \widehat{\otimes} \operatorname{Frac} k[[a]] \to \dual \widehat{\otimes} \operatorname{Frac} k[[a]], \\
	(f,g) \mapsto (1-a^{-1}t, 1) \cdot (f,g)
\end{gather*}
is an isomorphism.  We can conclude that $M=R+R \cdot (t^{2},0)$ satisfies Hypothesis~\ref{Hypothesis}. 

We deform the module $R \cdot (1,0) + R \cdot (t,1)$ to $R$ using the map 
\begin{gather*}
	\phi_{a} \colon R + R \cdot (t,0) \otimes k[[a]] \to k[[a]], \\
	f  \mapsto f_0 - a  f_1.
\end{gather*}
This map is a surjective homomorphism, and the kernel $\ker \phi_{a}$ is a deformation of $R \cdot (1,0) + R \cdot (t,1)$ to $R$.
Indeed, an isomorphism of the special fiber is given by
\begin{gather*}
	R \cdot (1,0) + R \cdot (t,1) \to M_{a} \otimes k[[a]]/(a), \\
	(f,g) \mapsto (t,t) \cdot (f,g).
\end{gather*}
and an isomorphism of the completed generic fiber is given by
\begin{gather*}
	M_{a} \widehat{\otimes} \operatorname{Frac} k[[a]] \to R \widehat{\otimes} \operatorname{Frac} k[[a]], \\
	(f,g) \mapsto (a^2-at+t^2, a^2) \cdot (f,g).
\end{gather*}
This proves that $M=R \cdot (1,0) + R \cdot (t,1)$ satisfies Hypothesis~\ref{Hypothesis}.

\subsubsection*{The $E_8(1)$-Singularity}  We need to deform the modules $R + R \cdot t^{4}$ and 
$R + R \cdot t$.  We begin with $R + R \cdot t^{4}$.

The module $R + R \cdot t^{4}$  deforms to the  dualizing module $\dual=R + R \cdot t^{2}$.
One such deformation is given by the kernel $M_{a} := \ker \phi_{a}$ of the surjective homomorphism
\begin{gather*}
	\phi_{a} \colon R + R \cdot t^{2} + R \cdot t^{4} \otimes k[[a]] \to k[[a]], \\
	f \mapsto f_{2} - a f_{4}.
\end{gather*}
To see this is a suitable deformation, observe that the special fiber of $M_{a}$ is equal to $R+R \cdot t^{4}$, and the map
\begin{gather*}
	M_{a} \widehat{\otimes} \operatorname{Frac} k[[a]] \to R + R \cdot t^{2} \widehat{\otimes} \operatorname{Frac} k[[a]], \\
	f \mapsto (a-t^{2}) \cdot f
\end{gather*}
is an isomorphism.  

Finally, we deform $R + R \cdot t$ to $R$.  Define $M_{a}$ to be the kernel of 
\begin{gather*}
	\phi_{a} \colon R \otimes k[[a]] \to k[\epsilon_1,\epsilon_2]/(\epsilon_1, \epsilon_2)^2 \otimes k[[a]], \\
	f \mapsto f_0 + \epsilon_1 (f_3+a f_5) + \epsilon_2 (f_7).
\end{gather*}
This map is surjective.  Furthermore, if $M_{a} := \ker \phi_{a}$, then the maps
\begin{gather*}
	R + R \cdot t \to M_{a} \otimes k[[a]]/(a), \\
	f \mapsto t^{5} \cdot f
\end{gather*}
and
\begin{gather*}
	M_{a} \widehat{\otimes} \operatorname{Frac} k[[a]] \to R \widehat{\otimes} \operatorname{Frac} k[[a]], \\
	f \mapsto (a t^{-3} - t^{-1}+a^{-1} t) \cdot f
\end{gather*}
are isomorphisms.  This shows that $M=R + R \cdot t$ satisfies Hypothesis~\ref{Hypothesis}.  Because $R + R \cdot t$
was the last module that we needed to deform, the proof is now complete.  \end{proof}

\section{Tables} \label{Sec: Tables}
Here we list the non-planar curve singularities of finite representation type (Table~\ref{Table: CurveSing}) and the rank $1$, torsion-free modules over the  ring of such a  singularity (Table~\ref{Table: BigTable}).  Both classification results are derived from \cite{greuel85}, where Greuel and Kn{\"o}rrer enumerate the maximal CM modules over the ring of an ADE singularity  \cite[pp.~423--425]{greuel85}.  

Their work can be used to classify the curve singularities of finite representation type as follows.  The main result of their paper states that a curve singularity is of finite representation type if and only if it dominates an ADE singularity \cite[Satz~1]{greuel85}.  As a consequence, we can conclude that the rings of  singularities of finite representation type are the algebras of the form $\End(M)$ for $M$ a rank $1$, torsion-free sheaf over an ADE singularity.  Indeed, if $R \subset R'$ corresponds to a dominance relation, then $R'$ considered as a $R$-module satisfies $\End_{R}(R')=R'$.  

The singularities that arise in this manner are listed in Table~\ref{Table: CurveSing}.  The first column (``\texttt{Singularity}") is the name of the singularity.  The second column (``\texttt{Parameterization}") presents the ring of the singularity as the subring of a product of power series rings topologically generated by an explicit set of elements.   The names we use for singularities are the names  used in the literature on simple curve singularities (e.g.~\cite{fruhbis}).  Each singularity arises as a partial desingularization of a DE-singularity.  The singularity $A_{n} \vee L$ is a partial desingularization of the $D_{n+3}$-singularity.  Indeed, the local ring of the $A_{n} \vee L$-singularity is the extension of the local ring of the $D_{n+3}$-singularity that is generated by either $(0,t)$ or  $(0,0,t)$ (depending on the parity of $n$).  Similarly, the local ring of $E_{k}(1)$ is generated over the local ring of the $E_{k}$-singularity by $t^5$ for $k=6$, $(t^{4},0)$ for $k=7$, and $t^7$ for $k=8$. 

The rank $1$, torsion-free modules over a ring from Table~\ref{Table: CurveSing} are listed in Table~\ref{Table: BigTable}.  Again, this has been derived from \cite{greuel85}.  Given an extension $R \subset R'$ contained in $\operatorname{Frac} R$ with $R$ the ring of an ADE singularity, the rank $1$, torsion-free $R'$-modules are exactly the rank $1$, torsion-free $R$-modules $M$ with the property that the inclusion $R \subset \End(M)$ extends to an inclusion $R \subset \End(M)$. 

Table~\ref{Table: BigTable} was generated by checking this condition for the modules listed in \cite{greuel85} and should be read as follows.  The first column (``\texttt{Singularity}") lists a curve singularity from Table~\ref{Table: CurveSing}.  The module in the second column (``\texttt{Greuel--Kn{\"o}rrer module}") is a module over the ring of the singularity in the first column.  The module is presented as a submodule of the total ring of fractions $\operatorname{Frac}(R)$ generated by an explicit set of elements.  Some of the modules we list do not appear in \cite[pp.~423--425]{greuel85} because Greuel and Kn{\"o}rrer only list the indecomposable modules.  For an indecomposable module, the presentation in Table~\ref{Table: BigTable} is chosen to coincide with the presentation from \cite{greuel85}.  The third column (``\texttt{Endomorphism Ring}") lists the endomorphism ring of the module.  This ring is always a product of rings of ADE singularities and the singularities appearing in Table~\ref{Table: CurveSing}.  We write ``$\cup$" to indicate that the endomorphism ring is a product (or, geometrically, a disjoint union) of the listed singularities.  Finally, every listed $R$-module $M$ is either isomorphic to the ring $S := \End(M)$ (considered as a $R$-module) or the dualizing module $\dual_{S}$ of that ring.  In the fourth column (``\texttt{Isomorphic to a ring?}"), we write ``\texttt{Yes}" if $M$ is isomorphic to its endomorphism ring, and we leave the entry blank otherwise.  We leave the entry in the final column (``\texttt{Isomorphic to a dualizing module?}") blank if $M$ is isomorphic to its endomorphism ring, and otherwise we write ``\texttt{Yes}" as $M$ is then isomorphic to $\dual_{S}$.  Note that when $S$ is Gorenstein and $M$ is isomorphic to $S$, we leave the last entry, titled ``\texttt{Isomorphic to a dualizing module?}", blank even though $M$ is isomorphic to $\dual_{S}=S$.

\begin{table}[hp]
\caption{Non-planar singularities of finite representation type}  \label{Table: CurveSing}
{
\begin{center}
\begin{tabular}{ l l }
	\texttt{Singularity} 					& \texttt{Parameterization}	 \\
	\hline 
	$A_n \vee L$, $n \ge 2$ and even 	& $k[[(t^{n+1}, 0), (t^{2}, 0), (0,t)]]$  \\ 
	$A_n \vee L$, $n \ge 1$ and odd 	& $k[[(t^{(n+1)/2}, -t^{(n+1)/2}, 0), (t, t, 0), (0, 0, t)]]$ \\
	$E_{6}(1)$ 					& $k[[t^{3}, t^{4}, t^{5}]]$ \\
	$E_{7}(1)$ 					& $k[[(t^{2}, t), (t^{3},0), (t^{4}, 0)]]$ \\ 
	$E_{8}(1)$ 					& $k[[t^{3}, t^{5}, t^{7}]]$ \\ 
\end{tabular}
\end{center}
}
\end{table}

\newpage

\begin{table}[hp]
{\tiny
\begin{center}
\caption{Rank $1$, torsion-free modules} \label{Table: BigTable}
\begin{tabular}{l l l p{2.5cm} p{2.5cm}}
	\texttt{Singularity}			& \texttt{Greuel--Kn\"{o}rrer Module}			& \texttt{Endomorphism Ring}	& \texttt{Isomorphic to a ring?}	 & \texttt{Isomorphic to a dualizing module?} \\
	\hline \hline
	
	$A_{n} \vee L$, $n$ even		& $R + (t^{n+1},0) \cdot R$				& $A_{n} \vee L$	& Yes 					& \\
							& $R + (t^{n-1},0) \cdot R$				& $A_{n-2} \vee L$	& Yes 					& \\
							& \dots								& \dots			& \dots 					& \dots \\
							& $R + (t^{3},0) \cdot R$					& $A_{2} \vee L$	& Yes					&  \\
							& $R + (t,0) \cdot R$						& $A_{1}$			& Yes 					& \\
	\hline
							& $R \cdot (1,0) + R \cdot (t^{n-1},1)$		& $A_{n} \vee L$	& 						& Yes \\
							& $R \cdot (1,0) + R \cdot (t^{n-3},1)$		& $A_{n-2} \vee L$	& 						& Yes \\
							& \dots								& \dots			& \dots					& \dots \\
							& $R \cdot (1,0) + R \cdot (t^{3},1)$			& $A_{2} \vee L$	& 						& Yes \\
	\hline
							& $R + R \cdot (1,0)$					& $A_{n}\cup\text{sm}$	& Yes				&  \\
							& $R + R \cdot (1,0) + R \cdot (t^{n-1},0)$		& $A_{n-2}\cup\text{sm}$ & Yes				& \\
							& \dots								& \dots			& \dots 					& \dots \\
							& $R+ R \cdot (1,0) + R \cdot (t^{3},0)$		& $A_{2}\cup\text{sm}$	& Yes				&  \\
							& $R+ R \cdot (1,0) + R \cdot (t,0)$			& $\text{sm}\cup\text{sm}$	& Yes 				& \\
	\hline
	$A_{n} \vee L$, $n$ odd		& $R + R \cdot (t^{(n+1)/2},0,0) $			& $A_{n} \vee L$	& Yes					&  \\
							& $R + R \cdot (t^{(n-1)/2},0,0) $			& $A_{n-2} \vee L$	& Yes 					& \\
							& \dots								& \dots			& \dots 					& \dots \\
							& $R + R \cdot (t^{2}, 0, 0)$				& $A_{3} \vee L$	& Yes					&  \\
							& $R + R \cdot (t,0,0)$					& $A_{1} \vee L$	& Yes					&  \\
	\hline
							& $R \cdot (1,1,0) + R \cdot (t^{(n-1)/2},0,1)$	& $A_{n} \vee L$	& 						& Yes \\
							& $R \cdot (1,1,0) + R \cdot (t^{(n-3)/2},0,1) $	& $A_{n-2} \vee L$	& 						& Yes \\
							&  \dots								& \dots			& \dots					& \dots \\
							& $R \cdot (1,1,0) + R \cdot (t,0,1) $			& $A_{3} \vee L$	& 						& Yes \\
							& $R \cdot (1,1,0) + R \cdot (1,0,1) $			& $A_{1} \vee L$	& 						& Yes \\
	\hline
							& $R+ R \cdot (1,0,0)$					& $A_1 \cup \text{sm}$	& Yes				&  \\
	\hline
							& $R+ R \cdot (0,1,0)$					& $A_{1} \cup \text{sm}$	& Yes 				& \\
	\hline
							& $R + R \cdot (0,0,1)$					& $A_{n} \cup \text{sm}$	& Yes 				& \\
							& $R+ R \cdot (0,0,1)+ R \cdot(t^{(n-1)/2},0,0)$ 	& $A_{n-2} \cup \text{sm}$	& Yes 				& \\
							& \dots						& \dots			& \dots 							& \dots \\
							& $R + R \cdot (0,0,1)+ R \cdot (t^{2},0,0) $	& $A_{3} \cup \text{sm}$	& Yes				& \\
							& $R + R \cdot (0,0,1) + R \cdot (t,0,0)$		& $A_1 \cup \text{sm}$	& Yes 				& \\
							& $\widetilde{R}$							& $\text{sm} \cup \text{sm} \cup \text{sm}$	& Yes	& \\
	\hline

	$E_{6}(1)$				& $R + R \cdot t^{5}$						& $E_{6}(1)$					& Yes 		& \\
							& $R + R \cdot t^2$ 						& $A_{2}$						& Yes 		& \\
							& $R + R \cdot t + R \cdot t^2$				& $\text{sm}$					& Yes 		& \\
	\hline
							& $R + R \cdot t$						& $E_{6}(1)$					& 			& Yes \\
	\hline 
	$E_{7}(1)$				& $R + R \cdot(t^{4},0)$					& $E_{7}(1)$					& Yes 		& \\
							& $R + R \cdot (t^{2},0)$					&$A_{2} \vee L$				& Yes 		& \\
							& $R + R \cdot(t,0) + R\cdot(t^2,0)$ 			& $A_{1}$						& Yes 		& \\
	\hline
							& $R+R \cdot (t,0)$						& $E_{7}(1)$					& 			& Yes \\
							& $R \cdot(1,0) + R \cdot (t,1)$				& $A_{2} \vee L$				& 			& Yes \\
	\hline
							& $R + R \cdot(1,0)$						& $A_{2} \cup \text{sm}$			& Yes 		& \\
							& $R + R \cdot(1,0) + R \cdot(t,0)$			& $\text{sm} \cup \text{sm}$		& Yes		&  \\
	\hline
	$E_{8}(1)$				& $R+ R \cdot t^{7}$						& $E_{8}(1)$					& Yes 		& \\
							& $R + R \cdot t^{4}$						& $E_{6}(1)$					& Yes 		& \\
							& $R + R \cdot t^{2} + R \cdot t^{4}$			& $A_{2}$						& Yes 		& \\
							& $\widetilde{R}$							& $\text{sm}$					& Yes 	& \\
	\hline
							& $R + R \cdot t^{2}$						& $E_{8}(1)$					& 			& Yes \\
							& $R + R \cdot t$						& $E_{6}(1)$					& 			& Yes \\
\end{tabular}
\end{center}
}
\end{table}

\section*{Acknowledgements}
We would like to thank the anonymous referee, Daniel Erman, Eduardo Esteves, Steven Kleiman, Robert Lazarsfeld, Louisa McClintock, and Filippo Viviani for helpful comments concerning exposition.

\bibliography{bibl}

\end{document}